\documentclass[12pt]{article}

\newtheorem{conjecture}{Conjecture}

\usepackage{amscd}
\usepackage{epsfig}
\usepackage{url}
\usepackage{amssymb}
\usepackage{amsmath}
\usepackage{amsfonts}
\usepackage{graphicx}

\def\Dbar{\leavevmode\lower.6ex\hbox to 0pt{\hskip-.23ex \accent"16\hss}D}

\def\bZ{{\mbox{\bf Z}}}

\begin{document}

{\bf\LARGE
\begin{center}
Construction of symmetric Hadamard matrices of order $4v$ 
for $v=47,73,113$
\end{center}
}

\noindent
{\bf N. A. Balonin$^a$}, Dr. Sc., Tech., Professor, korbendfs@mail.ru \\
{\bf D. {\v{Z}}. {\Dbar}okovi{\'c}$^b$}, PhD, Distinguished Professor Emeritus, djokovic@uwaterloo.ca \\
{\bf D. A. Karbovskiy$^a$}, alra@inbox.ru \\
\noindent
${}^{a}$Saint-Petersburg State University of Aerospace Instrumentation, 
67, B. Morskaia St., 190000, Saint-Petersburg, Russian Federation \\
${}^{b}$University of Waterloo, Department of Pure Mathematics and Institute for Quantum Computing, Waterloo, Ontario, N2L 3G1, Canada \\

\begin{abstract}
We continue our systematic search for symmetric Hadamard matrices 
based on the so called propus construction. In a previous paper 
this search covered the orders $4v$ with odd $v\le41$. In this 
paper we cover the cases $v=43,45,47,49,51$.
The odd integers $v<120$ for which no symmetric Hadamard matrices of order $4v$ are known are the following:
$$47,59,65,67,73,81,89,93,101,103,107,109,113,119.$$
By using the propus construction, we found several symmetric Hadamard matrices of order $4v$ for $v=47,73,113$. 

\end{abstract}

\section{Introduction}

In this paper we continue the systematic investigation,  
begun in  \cite{BBDKM:2017-a}, of the propus construction of symmetric Hadamard matrices.

Let us recall that a {\em Hadamard matrix} is a $\{1,-1\}$-matrix $H$  of order $m$ whose rows are mutually orthogonal, i.e. $HH^T=mI_m$, where $I_m$ is the identity matrix of order 
$m$. We say that $H$ is {\em skew-Hadamard matrix} if also $H+H^T=2I_m$.
The famous {\em Hadamard conjecture} asserts that Hadamard matrices exist for all orders $m$ which are multiples of $4$. (They also exist for $m=1,2$.) Similar conjectures have been proposed for symmetric Hadamard matrices and skew-Hadamard matrices, see e.g.  \cite[V.1.4]{CK:2007}. The smallest orders $4v$ for which such matrices have not been constructed are 668 for Hadamard matrices, 276 for skew-Hadamard matrices, and 188 for symmetric Hadamard matrices. Let us also mention that symmetric Hadamard matrices of orders 116, 156, 172 have been constructed only very recently, see 
\cite{DDK:SpecMatC:2015,BBDKM:2017-a}.

Since the size of a Hadamard matrix or a skew or symmetric 
Hadamard matrix can always be doubled, while preserving its type, we are interested mostly in the case where these matrices have order $4v$ with $v$ odd. 

The propus construction is based on the so called {\em Propus array}
\begin{equation} \label{Propus-array}
H=\left[ \begin{array}{cccc}
-C_1 & C_2R & C_3R & C_4R \\
C_3R & RC_4 & C_1 & -RC_2 \\
C_2R & C_1 & -RC_4 & RC_3 \\
C_4R & -RC_3 & RC_2 & C_1
\end{array} \right].
\end{equation}
In this paper, except in section \ref{sec:Exceptional}, the matrices $C_i$ will be circulants of order $v$ and the matrix 
$R$ will be the back-circulant identity matrix of order $v$,
$$
\label{Matrix-R}
R=\left[ \begin{array}{ccccc}
0 & 0 & \cdots & 0 & 1 \\
0 & 0 &        & 1 & 0 \\
\vdots &       &   &   \\
0 & 1 &        & 0 & 0 \\
1 & 0 &        & 0 & 0
\end{array} \right].
$$
The matrix $H$ will be a Hadamard matrix if 
\begin{equation} \label{eq:uslov-C}
\sum_i C_i C_i^T = 4vI_{4v}.
\end{equation}
(Superscript T denotes transposition of matrices.)
If also $C_1^T=C_1$ and $C_2=C_3$ then $H$ will be a symmetric 
Hadamard matrix. 

To construct the circulants $C_i$ satisfying the above conditions we use the cyclic propus difference families 
$(A_1,A_2,A_3,A_4)$ with parameters 
$(v;k_1,k_2,k_3,k_4;\lambda)$ such that $A_2=A_3$ and at least one of the base blocks $A_1,A_4$ is symmetric. The parameters must satisfy the three equations
\begin{eqnarray} \label{eq:osn}
&& \sum_{i=1}^4 k_i(k_i-1) = \lambda(v-1), \\
&& \sum_{i=1}^4 k_i = \lambda+v, \label{eq:gs} \\
&& k_2 = k_3. \label{eq:23}
\end{eqnarray}
We refer to such parameter sets as the {\em propus parameter sets}.

For the definitions of the terms that we use here
and the facts we mention below, we refer the reader to \cite{BBDKM:2017-a}. Without any loss of generality, we impose the following additional restrictions: 
\begin{equation} \label{eq:add}
v/2 \ge k_1,k_2; \quad k_1 \ge k_4.
\end{equation}
For convenience we say that the propus parameter sets satisfying 
these additional conditions are {\em normalized}.

For a given odd $v$ there exist at least one normalized propus  parameter set, see \cite[Theorem 1]{BBDKM:2017-a}. However, there exist even $v$ for which this is not true, see \cite[Theorem 2]{BBDKM:2017-a}.

It is conjectured in \cite{BBDKM:2017-a} that for each odd $v$ 
there exists at least one propus difference family in the 
cyclic group $\bZ_v$ of integers modulo $v$. But this may 
fail if we specify not only $v$ (odd) but also the  
parameters $k_1,k_2=k_3,k_4$.  
Our computations suggest that these exceptional propus 
parameter sets must have all $k_i$ equal to each other. For instance, there is no cyclic propus difference family having the parameters $(25;10,10,10,10;15)$. (This is also true for the propus difference families over the elementary abelian group 
$\bZ_5\times\bZ_5$.)

One of the authors developed a computer program to search 
for propus difference families. For the description of the 
algorithm used in the program we refer the reader to 
\cite{BBDKM:2017-a}. We used that program on 
PCs to construct many such families for odd (or even) $v$. 
The first version of the program was used in the range 
$v<43$. The second, improved version, was capable of finding
solutions for $v\le 51$. Some of the timings for these computations are given in section \ref{sec:app}.

In section \ref{sec:v=47} we give several examples of symmetric Hadamard matrices of new orders 188, 292, and 452. 

In section \ref{Parametri} we list the normalized propus parameter sets for odd $v\in\{43,45,\ldots,59\}$ and for each 
of them we indicate whether propus families with that parameter set exist and, if they do, which of the blocks A or D can 
be chosen to be symmetric. This list together with a similar list in \cite{BBDKM:2017-a} shows that there is a rich supply of 
propus type symmetric Hadamard matrices for orders $4v$ with
odd $v<50$. Sporadic examples are also known for $v=53,55,57$. The first undecided case is $v=59$.

In section \ref{sec:Exceptional} we focus on the case where $v=s^2$ is an odd square. We count the number of propus parameter sets $(v;x,y,y,z;\lambda)$ with $v=s^2$ by dropping the normalization condition $x\ge z$. This number, $N_s$, is also the number of positive odd integer solutions of a simple quadratic Diophantine equation, namely (\ref{eq:Dioph}). When 
$s$ is an odd prime then we conjecture that 
$N_s-s-1\in\{+1,-1\}$. We refer to the cases where $v$ is odd and all $k_i$ are equal as exceptional cases. They occur only when $v=s^2$. We also conjecture that every prime 
$s\equiv1 \pmod{4}$ can be written uniquely as 
$s=(a^2+b^2)/(a-b)$ where $a$ and $b$ are positive integers and $1<a\le (s-1)/2$. Moreover, the denominator $a-b$ is either a square or 2 times a square.

Finally in section \ref{sec:app}, for each of the normalized propus parameter sets with odd $v=43,45,\ldots,51$, but 
excluding the exceptional parameter set $(49;21,21,21,21;35)$, we list one or two examples of propus difference families.

\section{Symmetric Hadamard matrices of new orders}
\label{sec:v=47}

The smallest order $4v$ for which no symmetric Hadamard
matrix was known previously is $188=4\cdot47$.
There are four propus parameter sets 
$$
(47; 20, 22, 22, 18; 35), (47; 22, 20, 20, 19; 34), 
(47; 23, 19, 19, 21; 35),(47; 23, 22, 22, 17; 37)
$$
with $v=47$. In each case we constructed many such matrices, but here we record just two examples for each parameter set. In all four cases, $A$ is symmetric in the first and $D$ symmetric in the second example. As $B=C$ we omit the block $C$. 
The examples are separated by semicolons.

\begin{eqnarray*}
&& (47; 20, 22, 22, 18; 35) \\
&& [1,2,6,7,12,14,15,18,22,23,24,25,29,32,33,35,40,41,45,46],\\
&& [0,1,2,3,4,7,9,10,13,14,19,26,28,30,32,34,35,36,37,39,42,46],\\
&& [0,1,2,10,12,15,20,23,26,27,28,30,33,34,39,42,43,45];\\
&& [0,1,3,4,6,7,10,11,13,15,18,19,24,29,31,33,35,37,38,45],\\
&& [0,1,2,5,8,9,10,12,13,18,19,23,24,25,27,29,31,32,38,39,41,44],\\
&& [9,10,11,12,14,16,20,21,23,24,26,27,31,33,35,36,37,38];\\
&& \\
&& (47; 22, 20, 20, 19; 34) \\
&& [1,4,5,7,8,9,11,12,16,18,21,26,29,31,35,36,38,39,40,42,43,46],\\
&& [0,1,2,3,7,15,16,19,21,23,26,27,28,29,30,32,34,37,38,44],\\
&& [0,1,2,3,8,9,11,12,13,18,20,21,26,27,32,34,36,41,44];\\
&& [0,1,2,3,7,9,10,12,14,16,17,18,20,23,26,27,28,35,37,42,43,45],\\
&& [0,1,2,3,9,11,13,14,19,23,26,27,29,30,32,33,34,35,38,43],\\
&& [0,4,6,11,15,16,18,19,22,23,24,25,28,29,31,32,36,41,43];\\
&& \\
&& (47; 23, 19, 19, 21; 35) \\
&& [0,1,4,5,7,9,10,11,12,15,19,22,25,28,32,35,36,37,38,40,42,43,46],\\
&& [0,1,2,3,5,8,9,10,12,13,17,19,22,24,28,30,34,36,37],\\
&& [0,1,2,3,13,14,17,18,19,21,25,26,27,30,32,34,35,40,41,43,44];\\
&& [0,1,2,3,5,6,7,8,11,13,15,16,18,22,23,24,26,27,29,33,38,40,45],\\
&& [0,1,2,3,5,11,12,17,22,25,29,30,31,33,34,35,37,38,41],\\
&& [0,2,4,7,8,10,11,16,17,21,23,24,26,30,31,36,37,39,40,43,45];\\
&& \\
&& (47; 23, 22, 22, 17; 37)  \\
&& [0,1,4,11,12,13,15,16,18,19,21,22,25,26,28,29,31,32,34,35,36,43,46],\\
&& [0,1,2,3,4,5,8,10,13,18,19,21,23,25,27,29,30,36,39,41,42,43],\\
&& [0,1,2,5,6,10,11,12,15,21,25,26,33,38,40,41,45];\\
&& [0,1,2,3,4,7,8,10,11,12,13,16,19,21,23,25,26,29,31,33,34,35,45],\\
&& [0,1,2,3,4,8,9,12,13,14,17,18,19,20,26,27,29,31,34,37,40,44],\\
&& [0,2,6,13,15,18,20,21,22,25,26,27,29,32,34,41,45].\\
&& \\
\end{eqnarray*}

Let us give a concrete example. We choose the first parameter set above, $(47;20,22,22,18;35)$, and its first propus difference family, namely: 

\begin{eqnarray*}
A &=& [1,2,6,7,12,14,15,18,22,23,24,25,29,32,33,35,40,41,45,46],\\
B=C &=& [0,1,2,3,4,7,9,10,13,14,19,26,28,30,32,34,35,36,37,39,42,46],\\
D &=& [0,1,2,10,12,15,20,23,26,27,28,30,33,34,39,42,43,45].
\end{eqnarray*}
The binary $\{+1,-1\}$-sequences $a,b=c,d$ associated with the base blocks $A,B=C,D$ are:
\begin{eqnarray*}
a &=& [1,-1,-1,1,1,1,-1,-1,1,1,1,1,-1,1,-1,-1,1,1,-1,1,1,1,-1, \\
&& -1,-1,-1,1,1,1,-1,1,1,-1,-1,1,-1,1,1,1,1,-1,-1,1,1,1,-1,-1], \\
b=c &=& [-1,-1,-1,-1,-1,1,1,-1,1,-1,-1,1,1,-1,-1,1,1,1,1,-1,1,1,1, \\
&& 1,1,1,-1,1,-1,1,-1,1,-1,1,-1,-1,-1,-1,1,-1,1,1,-1,1,1,1,-1],\\
d &=& [-1,-1,-1,1,1,1,1,1,1,1,-1,1,-1,1,1,-1,1,1,1,1,-1,1,1,-1, \\
&& 1,1,-1,-1,-1,1,-1,1,1,-1,-1,1,1,1,1,-1,1,1,-1,-1,1,-1,1].
\end{eqnarray*}
The circulant matrices $C_1,C_2=C_3,C_4$ whose first rows are the sequences $a,b=c,d$ are depicted on Figure 1. Our convention is that a white square represents $+1$ and the black square 
$-1$. Note that $C_1$ is symmetric.

\begin{figure}
\includegraphics[width=\linewidth]{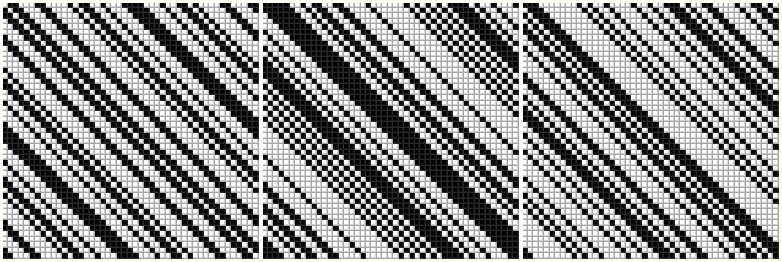}
\caption{The circulants $C_1,C_2=C_3,C_4$}
\label{fig:Fig1}
\end{figure}
By plugging these circulants into the propus array 
(\ref{Propus-array}) we obtain the symmetric Hadamard matrix  shown in Figure 2.

\begin{figure}
\begin{center}
\includegraphics[width=12cm]{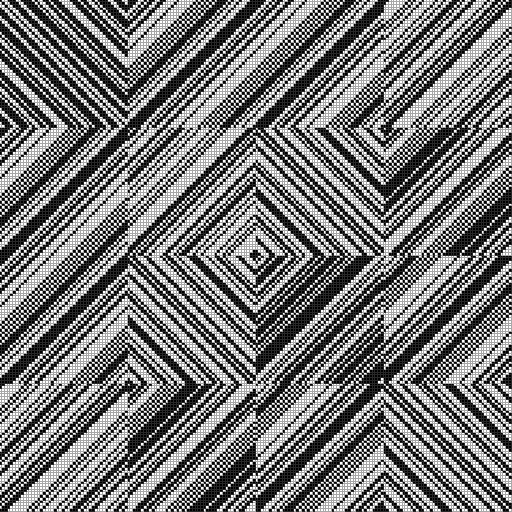}
\end{center}
\caption{Symmetric Hadamard matrix of order 188}
\label{fig:Fig2}
\end{figure}

Next we give the symmetric Hadamard matrices of order $4v$ 
where $v=73,113$. No symmetric Hadamard matrices of these orders 
were known previously.

We build the four base blocks $A,B=C,D$ as a union of orbits of a subgroup $H$ of $\bZ^*_v$ acting on the finite field $\bZ_v$. We choose $H=\{1,8,64\}$ for $v=73$ and $H=\{1,16,28,30,49,106,109\}$ for $v=113$.

For $v=73$ we use the parameter set $(73;36,36,36,28;63)$. 
The base blocks of the two propus difference families are:

\begin{eqnarray*}
A &=& \bigcup_{i\in I} iH,\quad 
I=\{1,2,3,4,9,11,18,21,26,27,36,43\} \\
B=C &=& \bigcup_{j\in J} jH,\quad 
J=\{2,4,5,6,9,12,14,17,27,34,35,36\} \\
D &=& \{0\} \cup \bigcup_{k\in K} kH,\quad 
K=\{1,2,3,6,7,9,18,42,43\}; \\
&& \\
A &=& \bigcup_{i\in I} iH,\quad 
I=\{1,2,3,4,6,9,12,18,25,27,35,36\} \\
B=C &=& \bigcup_{j\in J} jH,\quad 
J=\{2,5,7,9,13,17,25,26,33,35,36,42\} \\
D &=& \{0\} \cup \bigcup_{k\in K} kH,\quad 
K=\{4,6,13,18,27,34,35,36,42\}. \\
\end{eqnarray*}

For $v=113$ we use the parameter set $(113;56,49,49,56;97)$.
The base blocks of the four propus difference families are:

\begin{eqnarray*}
A &=& \bigcup_{i\in I} iH,\quad I=\{1,4,5,6,13,17,18,20\} \\
B=C &=& \bigcup_{j\in J} jH,\quad J=\{1,5,9,11,12,17,39\} \\
D &=& \bigcup_{k\in K} kH,\quad K=\{2,3,5,10,11,12,18,20\}; \\
&& \\
A &=& \bigcup_{i\in I} iH,\quad I=\{1,4,5,6,13,17,18,20\} \\
B=C &=& \bigcup_{j\in J} jH,\quad J=\{1,2,4,11,12,13,17\} \\
D &=& \bigcup_{k\in K} kH,\quad K=\{1,2,3,5,11,12,18,20\}; \\
&& \\
\end{eqnarray*}

\begin{eqnarray*}
A &=& \bigcup_{i\in I} iH,\quad I=\{1,4,5,6,13,17,18,20\} \\
B=C &=& \bigcup_{j\in J} jH,\quad J=\{1,2,4,11,12,13,17\} \\
D &=& \bigcup_{k\in K} kH,\quad K=\{3,4,5,8,9,12,13,20\}; \\
&& \\
A &=& \bigcup_{i\in I} iH,\quad I=\{1,3,4,10,12,13,18,39\} \\
B=C &=& \bigcup_{j\in J} jH,\quad J=\{2,5,9,10,17,20,39\} \\
D &=& \bigcup_{k\in K} kH,\quad K=\{2,3,9,11,12,17,20,39\}. \\
&& \\
\end{eqnarray*}

The first three families share the same block $A$,
and the second and third family differ only in block $D$. 
In spite of that, the four families are pairwise 
nonequivalent.

\section{Normalized propus parameter sets}
\label{Parametri}

We list here all normalized propus parameter sets 
$(v;x,y,y,z;\lambda)$ for odd $v=43,45,\ldots,59$. The cyclic propus families consisting of four base blocks 
$A,B,C,D\subseteq\bZ_v$ having sizes $x,y,y,z$, respectively, and such that $B=C$ and $A$ or $D$ is symmetric give symmetric Hadamard matrices of order $4v$. (If only $D$ is symmetric we have to switch $A$ and $D$ before plugging the blocks into the propus array.) If $x=z\ne y$ then the parameter set 
$(v;y,x,x,y;\lambda)$ is also normalized and is included in our list. In the former case the two base blocks of size $y$ have to be equal, while in the latter case the base blocks of size $x$ 
have to be equal. 

The four base blocks, subsets of $\bZ_v$,  are denoted by $A,B,C,D$. We require all propus difference families to have 
$B=C$. If we know such a family exists with symmetric block $A$, we indicate this by writing the symbol $A$ after the parameter set, and similarly for the symbol $D$. If we know that there exists a propus family with both $A$ and $D$ symmetric, then we write the symbol $AD$. Finally, the question mark means that the existence of a cyclic propus difference familiy remains undecided.

The symbol $T$ indicates that the parameter set belongs to the Turyn series of Williamson matrices. Since in that case all four base blocks are symmetric, the symbol $T$ implies $AD$. Further, the symbol $X$ indicates that the parameter set belongs to another infinite series (see 
\cite[Theorem 5]{DDK:SpecMatC:2015}) which is based on the 
paper \cite{XXSW} of Xia, Xia, Seberry, and Wu. In our list 
below the symbol $X$ implies $D$. More precisely, for a difference family $A,B,C,D$ in the $X$-series two blocks are equal, say $B=C$, 
and one of the remaining blocks is skew, block $A$ in our list, 
and the last one is symmetric, block $D$. 

For odd $v$ in the range $43,45,\ldots,51$ there is only one propus parameter set, $(49;21,21,21,21;36)$, for which we failed to find a cyclic propus difference family. 
(We believe that such family does not exist.)

\begin{center}
Normalized propus paramater sets with $v$ odd, $43\le v\le 59$
\end{center}

$$
\begin{array}{llll}
(43;18,21,21,16;33) & A,D    & (43;19,18,18,18;30) & A,D \\ 
(43;21,17,17,20;32) & A,D    & (43;21,19,19,16;32) & A,D \\ 
(43;21,21,21,15;35) & A,D    & (45;18,21,21,18;33) & A,D \\ 
(45;19,20,20,18;32) & AD,T   & (45;21,18,18,21;33) & A,D \\ 
(45;21,20,20,17;33) & A,D    & (45;21,22,22,16;36) & A,D \\ 
(45;22,19,19,18;33) & A,D,X  & (47;20,22,22,18;35) & A,D \\ 
(47;22,20,20,19;34) & A,D    & (47;23,19,19,21;35) & A,D \\ 
(47;23,22,22,17;37) & A,D    & (49;21,21,21,21;35) & ?   \\ 
(49;22,22,22,19;36) & A,D    & (49;22,24,24,18;39) & A,D \\ 
(49;23,20,20,22;36) & AD,T   & (49;23,23,23,18;38) & A,D \\ 
(51;21,25,25,20;40) & AD,T   & (51;23,22,22,21;37) & A,D \\ 
(53;22,24,24,22;39) & ?      & (53;24,22,22,24;39) & ?  \\
(53;24,25,25,20;41) & ?      & (53;26,22,22,23;40) & D,X \\
(55;23,26,26,22;42) & AD,T   & (55;24,25,25,22;41) & ?  \\
(55;24,27,27,21;44) & ?      & (55;26,23,23,24;41) & ?  \\
(55;27,24,24,22;42) & ?      & (55;27,25,25,21;43) & ?  \\
(57;25,25,25,24;42) & ?      & (57;27,25,25,23;43) & ?  \\
(57;27,26,26,22;44) & ?      & (57;28,28,28,21;48) & D,X \\
(59;26,28,28,23;46) & ?      & (59;27,25,25,26;44) & ?  \\ 
(59;28,29,29,22;49) & ?      & & \\ 
\end{array}
$$

In order to justify the claims made in this list, we give in section \ref{sec:app} examples of the propus difference families having the required properties. (For $v=47$ the examples are listed in section \ref{sec:v=47}.)

\section{Exceptional series of propus parameter sets}
\label{sec:Exceptional}

We say that a propus parameter set $(v;k_1,k_2,k_3,k_4:\lambda)$ is {\em exceptional} if $k_1=k_2=k_3=k_4$. The exceptional parameter sets are parametrized by just one integer $s>1$ and are given by the formula
\begin{equation} \label{def:exceptional}
\Pi_s=(~s^2;~\binom{s}{2},\binom{s}{2},\binom{s}{2},
\binom{s}{2};~s(s-2)~).
\end{equation}

There exists a cyclic propus difference family with parameter set $\Pi_3$. There exists also a propus difference 
family $(A,B,C,D)$ over the group $\bZ_3\times\bZ_3$ with 
the same parameter set and such that $A$ is symmetric and 
$B=C=D$. By using the finite field $\bZ_3[\alpha]$ where 
$\alpha^2=-1$, we can take 

\begin{equation} \label{3x3}
A=\{0,\alpha,-\alpha\},\quad 
B=C=D=\{\alpha,1-\alpha,\alpha-1\}.
\end{equation} 

For $s=5$, it is reported in \cite{BBDKM:2017-a} that there are no propus difference families in $\bZ_{25}$ having $\Pi_5$ as its parameter set. We performed another exhaustive search and found no such families in $\bZ_5\times\bZ_5$.

For $s=7$, our non-exhaustive searches found no cyclic propus difference families having the parameter set $\Pi_7$. 
However, we found a cyclic difference family with parameter set $\Pi_7$ and $B=C$ with neither $A$ nor $D$ symmetric:

\begin{eqnarray*}
A&=& [0,1,2,3,4,8,11,12,14,19,21,24,26,27,29,37,38,41,44,45,46], \\
B&=&C= [0,1,2,3,5,7,11,14,15,17,24,27,28,29,32,35,38,43,44,45,47],\\
D&=& [0,1,2,5,6,8,10,11,12,14,16,18,21,22,23,30,31,32,36,37,41].\\
\end{eqnarray*}

While computing the propus parameter sets $(v;x,y,y,z;\lambda)$ in the case when $v=s^2$ is an odd square, we observed an interesting feature. Namely, if in the definition of normalized propus difference sets we drop only the condition that 
$x\ge z$ and if $s$ is an odd prime then the number, $N_s$, of 
such parmeter sets is either $s$ or $s+2$.
It follows from the proof of \cite[Theorem 1]{BBDKM:2017-a} that $N_s$ is equal to the number of odd positive integer solutions of the Diophantine equation
\begin{equation} \label{eq:Dioph}
\xi^2+2\eta^2+\zeta^2=4s^2.
\end{equation}

After making additional computations, we decided to propose the 
following conjecture.

\begin{conjecture} \label{cnj:primes}
For any odd prime $s$, $N_s-s-1\in\{+1,-1\}$.
\end{conjecture} 

We have verified our conjecture for all odd primes less than 10000. There are 1228 such primes. For 606  of them we have $N_s=s$ and for the remaining 622 we have $N_s=s+2$. Thus the sequence $N_s-s-1$ is a $\{+1,-1\}$-sequence when $s$ runs through odd primes $<10000$. We have sketched the partial sums of this sequence on Figure 3.

If $s$ is a prime congruent to $1 \pmod{4}$, we observed that apart from $\Pi_s$ there is another normalized propus parameter set with $v=s^2$ and $k_2=k_3=\binom{s}{2}$. Let us denote this new parameter set by
\begin{equation} \label{def:exc-comp}
\Pi'_s=(~s^2;~\binom{s}{2}+\alpha,\binom{s}{2},\binom{s}{2},
\binom{s}{2}-\beta;~s(s-2)+\alpha-\beta~).
\end{equation}
The integers $\alpha$ and $\beta$ are positive and satisfy the quadratic Diophantine equation
\begin{equation} \label{eq:Dioph-comp}
\alpha^2+\beta^2=s(\alpha-\beta).
\end{equation}

We propose another conjecture.

\begin{conjecture} \label{cnj:primes1mod4}
For any odd prime $s\equiv1 \pmod{4}$ the Diophantine equation 
(\ref{eq:Dioph-comp}), in the unknowns $\alpha$ and $\beta$, has 
a unique solution $(a,b)$, where $a$ and $b$ are positive integers and $1<a\le (s-1)/2$. 
Moreover, $a-b$ is either a square or 2 times a square.
\end{conjecture} 

We have verified that this conjecture holds for $s<100000$.
If we drop the condition $1<a\le (s-1)/2$, than there exists one more solution, namely $(s-a,b)$.
Note also that the two solutions share the same $b$, and so 
the integer $b$ is uniquely determined by $s$.

\begin{figure}
\includegraphics[width=\linewidth]{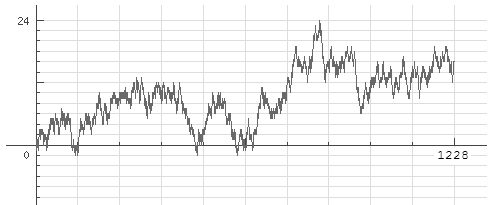}
\caption{Partial sums of the sequence $N_s-s-1$, $s$ odd prime}
\label{fig:Fig3}
\end{figure}

\newpage

\section{Appendix} \label{sec:app}

The cyclic propus difference families listed below, except 
some of the families that belong to one of the two infinite series $T$ and $X$, have been constructed by using a computer program written by one of the authors. The program was run on two PCs, each with a single 64-bit processor.
For $v=39$ it takes about 5 minutes to obtain a solution, 
about 20 minutes for $v=41$, about 1 hour for $v=43$, about 
3 or 4 hours for $v=45$, about 12 hours for $v=47$, about 
2 days for $v=49$, and 5 days for $v=51$.
In all families below the base block $B=C$, and to save space we omit the block $C$. The families are terminated by semicolons.

\begin{eqnarray*}
&& (43; 18, 21, 21, 16; 33) \\
&& [1,3,6,9,14,15,16,19,20,23,24,27,28,29,34,37,40,42],\\
&& [0,1,2,4,5,10,12,14,15,16,17,20,21,23,24,26,27,28,32,34,41],\\
&& [0,1,2,3,9,10,13,15,18,21,29,34,36,37,38,39];\\
&& [0,1,2,3,7,8,9,10,16,17,20,22,25,28,36,41],\\
&& [0,1,2,3,6,7,9,10,12,13,14,18,20,27,29,30,31,33,34,39,41],\\
&& [1,3,6,9,14,15,16,19,20,23,24,27,28,29,34,37,40,42];\\
&& \\
&& (43; 19, 18, 18, 18; 30) \\
&& [0,4,9,10,11,15,16,18,19,21,22,24,25,27,28,32,33,34,39],\\
&& [0,1,2,3,11,12,17,19,20,23,24,25,27,29,31,33,36,40],\\
&& [0,1,2,3,5,10,12,15,18,23,25,26,28,29,36,39,40,41];\\
&& [0,1,2,5,9,10,14,16,20,23,24,27,29,30,32,34,36,38,40],\\
&& [0,1,2,3,4,8,9,10,11,14,18,21,26,27,30,32,40,42],\\
&& [2,7,8,9,10,13,15,18,21,22,25,28,30,33,34,35,36,41];\\
&& \\
&& (43; 21, 17, 17, 20; 32) \\
&& [0,1,3,6,7,9,11,14,16,20,21,22,23,27,29,32,34,36,37,40,42],\\
&& [0,1,2,3,5,6,7,13,15,24,25,28,29,32,37,39,40],\\
&& [0,1,2,3,10,12,14,15,18,19,20,25,28,29,31,32,34,35,37,39];\\
&& [0,1,2,3,4,7,12,13,14,18,20,23,24,28,30,32,33,34,36,38,41],\\
&& [0,1,2,5,8,10,15,17,18,19,21,24,25,30,36,37,40],\\
&& [1,3,4,5,6,7,8,13,18,21,22,25,30,35,36,37,38,39,40,42];\\
&& \\
&& (43; 21, 19, 19, 16; 32)  \\
&& [0,1,6,11,12,13,16,17,19,20,21,22,23,24,26,27,30,31,32,37,42],\\
&& [0,1,2,6,8,9,12,15,17,20,22,23,24,26,27,35,36,39,41],\\
&& [0,1,2,6,8,9,11,15,16,18,20,24,28,29,31,41];\\
&& [0,1,2,3,4,7,11,13,15,17,19,20,22,32,33,34,35,37,39,40,42],\\
&& [0,1,2,4,5,8,10,11,13,16,17,20,23,24,25,27,34,38,39],\\
&& [2,3,4,10,12,14,15,20,23,28,29,31,33,39,40,41];\\
\end{eqnarray*}

\begin{eqnarray*}
&& (43;21,21,21,15;35) \\
&& [0,1,2,3,4,8,9,12,14,19,22,23,26,28,29,31,32,34,38,39,41],\\
&& [1,4,6,9,10,11,13,14,15,16,17,21,23,24,25,31,35,36,38,40,41],\\
&& [0,7,9,13,14,15,17,18,25,26,28,29,30,34,36]; \\
&& [0,1,2,3,6,7,9,11,13,15,16,17,21,22,25,26,29,33,38,39,41],\\
&& [0,1,2,3,4,5,6,7,8,13,16,17,20,22,25,27,29,34,35,37,40],\\
&& [0,5,6,12,13,14,16,20,23,27,29,30,31,37,38];\\
\end{eqnarray*}

The last example consists of a D-optimal design (blocks $A$ and $D$) and two copies of the Paley difference set in $\bZ_{43}$ 
(blocks $B=C$). It is taken from the paper \cite{DDK:SpecMatC:2015}.

\begin{eqnarray*}
&& (45; 18, 21, 21, 18; 33) \\
&& [4,7,8,9,10,11,16,19,20,25,26,29,34,35,36,37,38,41],\\
&& [0,1,2,3,5,7,8,12,13,16,19,22,23,27,32,34,36,39,40,42,44],\\
&& [0,1,2,3,10,12,13,15,17,19,24,25,32,34,37,38,39,41];\\
&& \\
&& (45; 19, 20, 20, 18; 32) \\
&& [0,1,6,12,13,14,16,17,20,22,23,25,28,29,31,32,33,39,44],\\
&& [1,3,7,8,10,11,12,13,17,20,25,28,32,33,34,35,37,38,42,44],\\
&& [1,6,12,13,14,16,17,20,22,23,25,28,29,31,32,33,39,44];\\
&& \\
&& (45; 21, 18, 18, 21; 33) \\
&& [0,4,6,11,12,13,14,17,18,20,22,23,25,27,28,31,32,33,34,39,41],\\
&& [0,1,2,5,6,8,9,11,13,20,21,23,31,32,34,37,38,41],\\
&& [0,1,2,3,6,10,12,17,18,19,21,22,23,25,33,35,37,38,40,41,42];\\
&& \\
&& (45; 21, 20, 20, 17; 33)  \\
&& [0,2,3,4,5,10,14,15,17,19,22,23,26,28,30,31,35,40,41,42,43],\\
&& [0,1,2,3,4,8,12,13,15,22,23,26,28,30,31,36,37,39,42,43],\\
&& [0,1,2,3,4,6,10,14,20,26,27,30,33,35,37,42,43];\\
&& [0,1,2,3,8,9,15,16,18,22,23,25,28,32,33,34,36,38,39,42,43],\\
&& [0,1,2,3,5,6,8,9,12,14,16,18,19,20,24,27,29,31,32,41],\\
&& [0,1,2,7,10,12,18,19,22,23,26,27,33,35,38,43,44];\\
\end{eqnarray*}

\begin{eqnarray*}
&& (45; 21, 22, 22, 16; 36)  \\
&& [0,3,4,6,7,9,11,12,13,14,18,27,31,32,33,34,36,38,39,41,42],\\
&& [0,1,2,3,4,6,9,11,14,15,16,19,20,26,28,30,34,35,36,38,42,43],\\
&& [0,1,2,3,10,11,13,17,22,23,24,27,30,34,39,42];\\
&& [0,1,2,4,6,7,9,11,12,17,21,24,25,27,28,30,32,33,34,39,43],\\
&& [0,1,2,4,5,9,10,13,14,15,18,20,21,26,28,29,35,36,38,40,42,43],\\
&& [3,9,13,15,16,17,18,19,26,27,28,29,30,32,36,42];\\
&& \\
&& (45; 22, 19, 19, 18; 33) \\
&& [2,3,4,7,8,10,12,13,14,17,20,25,28,31,32,33,35,37,38,41,42,43],\\
&& [0,1,2,5,9,11,14,18,19,20,22,24,26,27,30,31,32,33,34],\\
&& [0,1,2,7,9,10,13,16,17,19,24,27,33,35,36,38,40,43];\\
&& [0,1,2,3,7,10,11,15,16,18,19,20,25,28,30,31,35,36,37,40,42,43],\\
&& [0,1,2,4,6,12,19,20,21,24,25,29,31,32,33,35,40,42,43],\\
&& [1,3,4,5,8,10,11,18,21,24,27,34,35,37,40,41,42,44];\\
&& \\
&& (49; 22, 22, 22, 19; 36) \\
&& [1,3,5,8,9,11,12,15,16,18,19,30,31,33,34,37,38,40,41,44,46,48],\\
&& [0,1,2,3,4,5,6,9,14,15,18,25,27,30,32,33,35,37,38,42,43,44],\\
&& [0,1,2,5,6,10,12,14,18,19,21,27,32,34,35,36,40,43,45];\\
&& [0,1,2,3,4,7,10,14,15,18,19,24,26,30,31,32,33,35,37,40,41,47],\\
&& [0,1,2,3,4,6,10,11,14,16,17,23,24,31,34,36,38,39,41,42,43,47],\\
&& [0,1,3,6,8,12,18,21,22,23,26,27,28,31,37,41,43,46,48];\\
&& \\
&& (49; 22, 24, 24, 18; 39) \\
&& [2,3,6,8,9,17,19,20,21,22,24,25,27,28,29,30,32,40,41,43,46,47],\\
&& [0,1,2,3,7,8,9,16,19,21,23,25,26,28,29,32,34,36,37,38,40,
41,42,46],\\
&& [0,1,2,3,7,8,12,15,17,19,26,29,36,37,39,42,43,46];\\
&& [0,1,2,3,7,10,12,13,14,17,19,20,22,23,25,26,27,28,32,37,40,46],\\
&& [0,1,2,3,4,5,6,9,11,13,14,15,17,19,22,27,28,31,33,34,35,38,43,45],\\
&& [2,5,6,10,16,17,19,23,24,25,26,30,32,33,39,43,44,47];\\
\end{eqnarray*}

\begin{eqnarray*}
&& (49; 23, 20, 20, 22; 36)  \\
&& [0,1,2,4,6,8,15,16,17,20,21,23,26,28,29,32,33,34,41,43,45,47,
48], \\
&& [3,6,10,13,14,15,20,21,23,24,25,26,28,29,34,35,36,39,43,46],\\
&& [1,2,4,6,8,15,16,17,20,21,23,26,28,29,32,33,34,41,43,45,47,
48];\\
&& \\
&& (49; 23, 23, 23, 18; 38)  \\
&& [0,3,4,7,9,10,11,13,15,17,18,23,26,31,32,34,36,38,39,40,42,45,46],\\
&& [0,1,2,3,4,5,6,9,10,11,12,15,16,21,23,25,27,28,32,35,40,41,43],\\
&& [0,1,2,5,6,13,15,18,21,22,27,32,34,36,37,39,46,47];\\
&& [0,1,2,3,4,5,6,7,8,11,15,17,20,24,27,29,33,36,38,41,44,45,47],\\
&& [0,1,2,3,5,7,8,10,12,14,20,22,23,24,30,31,32,35,36,37,38,41,46],\\
&& [3,4,10,13,14,16,20,21,24,25,28,29,33,35,36,39,45,46];\\
&& \\
\end{eqnarray*}

\begin{eqnarray*}
&& (51;23,22,22,21;37) \\ 
&& [0,2,4,10,11,13,16,20,21,23,24,25,26,27,28,30,31,35,38,40,
41,47,49], \\
&& [0,1,2,3,4,5,9,13,15,20,22,23,27,31,32,33,38,39,41,44,47,48],\\
&& [0,1,2,3,6,8,9,12,13,14,17,19,22,31,34,36,37,38,40,44,49];\\
&& [0,1,2,3,4,5,10,12,13,14,15,19,21,22,28,30,34,37,39,41,42,
47,49],\\
&& [0,1,2,4,5,8,10,13,18,19,21,24,25,28,29,31,33,35,38,39,40,43],\\
&& [0,2,3,4,5,6,9,13,19,20,25,26,31,32,38,42,45,46,47,48,49];\\
&& (51;25,25,21,20;40) \\
&& [0,1,4,5,7,9,15,16,17,18,22,29,33,34,35,3642,44,46,47,50],\\
&& [0,1,2,4,5,7,8,11,15,16,21,23,25,26,28,30,35,36,40,43,44,46,
47,49,50],\\
&& [1,4,5,7,9,15,16,17,18,22,29,33,34,35,36,42,44,46,47,50];\\
&& (53;26,22,22,23;40) \\ 
&& [1,5,6,10,11,12,15,18,22,27,28,29,30,32,33,34,36,37,39,40,44,45,46,49,50,51], \\
&& [0,1,2,3,9,11,18,21,24,25,29,33,34,35,36,41,44,46,48,49,50,52],\\
&& [0,1,3,9,10,12,14,16,17,20,23,25,28,30,33,36,37,39,41,43,44,50,
52];\\
\end{eqnarray*}

\begin{eqnarray*}
&& (55;23,26,26,22;42)  \\
&& [0,6,7,10,11,15,17,18,19,21,24,26,29,31,34,36,37,38,40,44,45,48,49], \\
&& [1,2,4,8,14,16,17,18,19,23,24,25,27,28,30,31,32,36,37,38,39,41,
47,51,53,54],\\
&& [6,7,10,11,15,17,18,19,21,24,26,29,31,34,36,37,38,40,44,45,48,49];\\
&& \\
&& (57;28,28,28,21;48)  \\
&& [2,4,12,13,15,21,23,24,25,27,28,31,35,37,38,39,40,41,43,46,47,48,49,50,51,52,54,56],\\
&& [0,1,2,3,4,6,9,11,13,16,17,20,23,28,31,32,34,35,37,39,40,41,43,44,45,49,50,53],\\
&& [0,1,4,6,13,14,15,19,20,21,26,31,36,37,38,42,43,44,51,53,56];\\
&& \\
\end{eqnarray*}

\section{Acknowledgements}
The research of the first author leading to these results has received funding from the Ministry of Education and Science of the Russian Federation according to the project part of the state funding assignment No 2.2200.2017/4.6.
The second author wishes to acknowledge generous support by NSERC. His work was made possible by the facilities of the Shared Hierarchical Academic Research Computing Network (SHARCNET) and Compute/Calcul Canada.

\end{document}